\documentclass[12pt]{article}%
\usepackage{amsfonts}
\usepackage{amsmath,amssymb}
\usepackage[applemac]{inputenc}
\usepackage{amsmath,amssymb,fullpage}
\usepackage{color}
\usepackage{amsmath}
\usepackage{amssymb}
\usepackage{graphicx}
\usepackage[margin=1.5in]{geometry}%
\providecommand{\U}[1]{\protect\rule{.1in}{.1in}}

\newcommand{\dproof}{\noindent {Proof.} \quad}
\newcommand{\fproof}{\hfill $\square$ \bigskip}
\newtheorem{definition}{Definition}[section]
\newtheorem{notation}{Notation}[section]

\newtheorem{assumption}{Assumption}[section]
\newtheorem{theorem}[definition]{Theorem}

\newtheorem{remark}[definition]{ \it Remark}

\newtheorem{lemma}[definition]{Lemma}

\numberwithin{equation}{section}

\def\1B{\text{1\!\!I}}

\newcommand{\mc}{\mathcal}

\begin{document}

\date{\today}
\title{Stochastic Maximum Principle with Default}

\author{Khalida Bachir Cherif$^{1}$, Nacira Agram$^{2}$, Kristina Dahl$^{3}$ }

\maketitle

\begin{abstract}
In this paper, we derive sufficient and necessary maximum principles for a stochastic optimal control problem where the system state is given by a controlled stochastic differential equation driven by a Brownian motion and a pure jump martingale. We apply the maximum principles to solve a log-utility maximisation problem with default. 
\end{abstract}

\footnotetext[1]{{\small University D. Moulay Taher, Saida, Algeria. Email:
\texttt{khalida.bc@hotmail.fr.}}}

\footnotetext[2]{{\small Department of Mathematics, Linnaeus University
SE-351\,95 V\"axj\"o, Sweden. Email:\texttt{ nacira.agram@lnu.se.}}}

\footnotetext[3]{{\small Department of Mathematics, University of Oslo, O.
Box 1053 Blindern, N--0316 Oslo, Norway. Email:\texttt{ kristrd@math.uio.no.}
}}

\paragraph{Keywords:}

Stochastic maximum principle, Backward Stochastic Differential Equations with Default, Single Jump.\newline

\section{Introduction}

In this paper we study optimal control of stochastic systems with default.  By saying that the system is with default we mean that the system is driven by both the Brownian motion and the martingale $M$ associated with a default jump with an intensity process $\lambda$. We derive sufficient and necessary conditions of optimality. The associated adjoint process is a solution of a backward stochastic differential equations (BSDE) driven by both the Brownian motion and the pure jump martingale $M$. 

This type of equations has been studied by Dumitrescu et al. \cite{DumitrescuEtAl}. They prove existence and uniqueness as well as comparison theorems for these types of BSDEs. They also generalize the results to drivers including a singular process. If the driver is $\lambda$-linear, they find a representation of the solution of the associated BSDE in terms of a conditional expectation and an adjoint exponential semi-martingale. The framework of Dumitrescu et al. \cite{DumitrescuEtAl} is the same as that of our paper. However, in contrast to Dumitrescu et al. \cite{DumitrescuEtAl}, we consider a stochastic optimal control problem in this default framework, and derive maximum principles characterizing the optimal solution of this problem.

Several other papers have studied similar frameworks:\\
Kharroubi and Lim \cite{kha} consider BSDEs with random marked jumps as well as applications to default risk. They connect the BSDEs with random marked jumps to Brownian BSDEs by enlargement of filtrations and prove that the jump BSDEs have solutions if the Brownian BSDEs have solutions. Furthermore, they prove a uniqueness theorem for jump BSDEs via a comparison theorem with Brownian BSDEs. Though the framework of Kharroubi and Lim \cite{kha} is similar to that of our paper, we focus on the stochastic control problem instead of the properties of the BSDEs.\\
Lim and Quenez \cite{LimQuenez} analyze the exponential utility maximization problem for an incomplete market with a default time which causes a discontinuity in the stock price. They apply dynamic programming to characterize the value function as the maximal subsolution of a BSDE. Lim and Quenez \cite{LimQuenez2} consider a financial market with an asset exposed to a risk which can cause a jump in the asset price. They assume that the asset can be traded after the default time. In this context, they study the expected utility maximization of terminal wealth for several utility functions. They prove that the value function for the power utility function can be determined as the minimal solution of a BSDE.  Though the framework of Lim and Quenez \cite{LimQuenez} and \cite{LimQuenez2} is similar to ours, they consider the control problem for maximizing the utility of terminal wealth. In contrast, the objective function in the present paper consists of both a terminal time term and an integral term over the whole time period.\\
For information about stochastic control with default jumps, we refer to Pham \cite{P}.

For more on stochastic control for jump processes, see e.g. Cohen and Elliott \cite{CohenElliott} Chapter 21. Note that in Cohen and Elliott \cite{CohenElliott}, the control cannot affect the diffusion coefficient function. In the present paper, this is possible. 

This work is organised as follows: 
\begin{itemize}
\item
In Section 2 we give some preliminaries. 
\item
In Section 3 we study the stochastic maximum principle and we derive sufficient and necessary conditions for optimality. 
\item
Finally, we apply our results to solve logarithmic utility maximisation problem for a defaultable cash flow.
\end{itemize}

\section{Framework}

Let $(\Omega,\mathcal{G},P)$ be a complete probability space. We assume that
this space is equipped with a one-dimensional standard Brownian motion $W$ and
a single jump process $H_{t}=\mathbf{1}_{\tau\leq t},\ t\in\lbrack0,T],$ where the random variable $\tau$ is positive and may represent a default time in credit- or counterparty risk, or a death time in actuarial issues\footnote{If $\tau$ is a death time, the control is often stopped at this time. This complicates the problem, see e.g., Bouchard and Pham \cite{BouchardPham}, Choulli and Yansori \cite{ChoulliYansori} and Jeanblanc et al. \cite{JeanblancEtAl}.}. We assume that this default can appear at any time, i.e. $P(\tau\geq t)>0$ for any $t\geq0$. We denote by $\mathbb{G}:=(\mathcal{G}_{t})_{t\geq0}$ the complete natural filtration of $W$ and $H$. We assume that $W$ is a $\mathbb{G}$-Brownian motion.

We suppose that the increasing process $H$ admits a predictable compensator $\Lambda$. Moreover, the process $\Lambda$ is assumed to be absolutely continuous w.r.t. Lebesgue's measure, there exists a positive process  $\lambda$, called the intensity, such that $\Lambda_{t}=\int_{0}^{t}\lambda_s ds$ for each $t\geq0$.
\vspace{2mm}
The process $M$ defined as
\small
\begin{align*}
  M_t=H_t-\int_0^t \lambda_sds.  
\end{align*}
is a $\mathbb{G}$-martingale called the compensated martingale of $H$. If the intensity is $\mathbb{G}$-adapted, it vanishes after $tau$. This is important for the following.
\vspace{2mm}
We will state the so-called predictable representation theorem (PRT) (Theorem 3.12 in Aksamit and Jeanblanc \cite{aj}, reformulated to the current notation.)

\begin{theorem}Every ${\mathbb{G}}$-martingale $Y$ admits a
representation
\[
Y_{t}=Y_{0}+%
{\textstyle\int_{0}^{t}}
\varphi_{s}dW_{s}+%
{\textstyle\int_{0}^{t}}
\gamma_{s}dM_{s},
\]
where $M$ is the compensated martingale of $H$, and $\varphi=(\varphi_t)_{t\in [0,T]},\gamma=(\gamma_t)_{t\in [0,T]}$ are
${\mathbb{G}}$-predictable processes, such that the stochastic integrals are well defined.
\end{theorem}

Throughout this section, we introduce some basic spaces.

\begin{itemize}
\item $S^{2}$ is the subset of $%
\mathbb{R}
$-valued $\mathbb{G}$-adapted c{\`a}dl{\`a}g processes $\left(  Y_{t}\right)
_{t\in\left[  0,T\right]  }$, such that
\[
\left\Vert Y\right\Vert _{S^{2}}^{2}:=\mathbb{E[}\underset
{t\in\left[  0,T\right]  }{\sup}\left\vert Y_{t}\right\vert ^{2}]<\infty.
\]

\item $H^{2}$ is the subset of $%
\mathbb{R}
$-valued $\mathbb{G}$-predictable processes $\left(  Z_{t}\right)
_{t\in\left[  0,T\right]  },$ such that
\[
\left\Vert Z\right\Vert _{H^{2}}^{2}:=\mathbb{E[}%
{\textstyle\int_{0}^{T}}
\left\vert Z_{t}\right\vert ^{2}dt]<\infty.
\]

\item $H^{2}(\lambda)$ is the subset of $%
\mathbb{R}
$-valued $\mathbb{G}$-predictable processes $\left(  U_{t}\right)
_{t\in\left[  0,T\right]  },$ such that
\[
\left\Vert U\right\Vert _{H^{2}(\lambda)}^{2}:=\mathbb{E[}%
{\textstyle\int_{0}^{T}}
\lambda_{t}\left\vert U_{t}\right\vert ^{2}dt]<\infty.
\]

\end{itemize}
\section{Stochastic Maximum Principles}
In this section, we present two stochastic maximum principles which can be used to solve stochastic optimal control problems where the system state is determined by the controlled  with default.

Let $(u_t)_{t \geq 0}$ be a control process. We denote by $\mathcal{V}$ a given convex subset of $\mathbb{R}$ and we let $\mathcal{A}$ be the set of admissible controls. Assume that $\mathcal{A}$ is a given set of $\mathcal{V}$-valued, $\mathbb{G}$-predictable processes in $L^{2}(\Omega\times
\lbrack0,T])$.

\noindent Consider the controlled stochastic differential equation (SDE) with default under the filtration $\mathbb{G},$
\begin{equation}
\begin{array}
[c]{ll}%
dX_{t} & =b(t,X_{t},u_{t})dt+\sigma(t,X_{t},u_{t})dW_{t}+\gamma(t,X_{t^{-}}%
,u_{t})dM_{t};\quad
X_{0} =x_{0},
\end{array} \label{eq: default}%
\end{equation}

\noindent where the coefficient functions are as follows:
\begin{align*}
b &  :\Omega\times\lbrack0,T]\times\mathbb{R}\times\mathcal{V}\rightarrow
\mathbb{R},\\
\sigma &  :\Omega\times\lbrack0,T]\times\mathbb{R}\times\mathcal{V}%
\rightarrow\mathbb{R},\\
\gamma &  :\Omega\times\lbrack0,T]\times\mathbb{R}\times\mathcal{V}%
\rightarrow\mathbb{R},
\end{align*}
and the initial value $x_{0}\in%
\mathbb{R}
.$

\noindent Note that we often suppress the $\omega$ for ease of notation. So,
for instance, we write $b(t,X_{t},u_{t})$ instead of $b(\omega,t,X_{t}(\omega),u_{t}(\omega))$.

\noindent We make the following set of assumptions on these coefficient functions.

\begin{assumption}
\label{assumption: E}

\begin{enumerate}
\item[$(a)$] The functions $b(\omega,t,\cdot)$, $\sigma(\omega,t,\cdot)$ and
$\gamma(\omega,t,\cdot)$ are assumed to be bounded and $C^{1}$ for each fixed $\omega,t$ with bounded derivatives.

\item [$(b)$] The functions $b(\cdot,x,u)$ and $\sigma(\cdot,x,u)$ and $\gamma(\cdot,x,u)$ are $\mathbb{G}%
-$predictable, for each $(x,u)\in \mathbb{R} \times \mathcal{V}$.

\item [$(c)$] \emph{Lipschitz condition:} The functions $b,\sigma, \gamma$ are uniformly Lipschitz in the variable $x$ for each $u\in\mathcal{V}$, with the Lipschitz constant, $\psi >0$, independent of the
variables $t,\omega$.

\item [$(d)$] \emph{Linear growth:} The functions $b,\sigma,\gamma$ satisfy the linear
growth condition in the variable $x$,  for each $u\in\mathcal{V}$, with the linear growth constant
independent of the variables $t,\omega$.
\end{enumerate}
\end{assumption}
\begin{theorem}[Existence of unique solution to the SDE with default]
\label{prop} Under the assumptions (b-d), there exists a unique solution $X\in
S^{2}$ of SDE \eqref{eq: default}.
\end{theorem}
In fact, the proof of this theorem can be obtained as a consequence of a general result of Theorem 16.3.1 in Cohen and Elliott \cite{CohenElliott}.
Now that we know that there exists a unique solution to the controlled SDE with default \eqref{eq: default}, we can move on to study a stochastic optimal control problem with default.

The performance functional, which we would like to maximise over all
strategies $u\in\mathcal{A}$, is defined as
\[
J(u)=E\left[  \int_{0}^{T}h(t,X_{t},u_{t})dt+g(X_{T})\right]  .
\]
We assume that the functions
\begin{align*}
h &  :\Omega\times\lbrack0,T]\times\mathbb{R}\times\mathcal{V}\rightarrow
\mathbb{R},\\
g &  :\Omega\times\mathbb{R}
\rightarrow\mathbb{R},
\end{align*}
are $ \mathbb{G}$-predictable, $ \mathcal{G}_T$-measurable respectively, $C^1$  w.r.t. $x,u$ and admits bounded derivatives. Moreover, 
$$
E\left[  \int_{0}^{T}h^2(t,X_{t},u_{t})dt+g^2(X_{T})\right]<\infty.
$$
We would like to derive stochastic maximum principles for this problem. \\
The associated Hamiltonian functional is defined by
\begin{equation}
\label{eq: Hamiltonian}
\mathcal{H}(t,x,u,p,q,w):=h(t,x,u)+b(t,x,u)p+\sigma(t,x,u)q+\lambda_t\gamma(t,x,u)w,
\end{equation}
where $p,q,w$ are called the adjoint processes. 

\begin{notation}
\label{not: shorthand}
For ease of notation, we define the following shorthand for some given control $u$ with corresponding $X$
\[
b_t := b(t, X_t, u_t), \\[\smallskipamount]
\sigma_t := \sigma(t, X_t, u_t), \\[\smallskipamount]
\gamma_{t} := \gamma(t,X_{t^{-}}, u_t),
\]
\[
\hat{b}_t := b(t, \hat{X}_t, \hat{u}_t), \\[\smallskipamount]
\hat{\sigma}_t := \sigma(t, \hat{X}_t, \hat{u}_t), \\[\smallskipamount]
\hat{\gamma}_{t} := \gamma(t, \hat{X}_{t^{-}}, \hat{u}_t).\\
\]
We will use the same notations for the partial derivatives of the above coefficients.
\end{notation}

\noindent The adjoint processes $p,q,w$ are given as the solution of the adjoint BSDE%
\[
\begin{array}
[c]{ll}%
dp_{t} & =-\frac{\partial\mathcal{H}}{\partial x} (t, X_t, u_t, {p}_t, {q}_t, {w}_t)dt +q_{t}dW_{t}+w_{t}%
dM_{t};\quad
p_{T} =g^{\prime}(X_{T}).
\end{array}
\]
\noindent Using the definition of the Hamiltonian \eqref{eq: Hamiltonian}, the above adjoint BSDE can be rewritten as:%

\begin{equation}
\begin{array}
[c]{llll}%
dp_{t}=-\Big[\frac{\partial h_{t}}{\partial x}+\frac{\partial b_{t}}{\partial
x}p_{t}+\frac{\partial
\sigma_{t}}{\partial x}q_{t}+\lambda_t\frac{\partial\gamma_{t}}{\partial x}w_{t}\Big]dt+q_{t}dW_{t}+w_{t}dM_{t};\quad p_{T} =g^{\prime}(X_{T}).
\end{array}  \label{eq: adjoint}%
\end{equation}
Note that this adjoint equation is a linear BSDEs with default, which by our assumptions of the coefficients and Theorem 2.14 in Dumitrescu et al. \cite{DumitrescuEtAl}, has the following explicit solution

\[
 p_t= E\Big[\Gamma_{t,T}g^{\prime}(X_{T}) + \int_t^{T} \Gamma_{t,s} \frac{\partial h_{s}}{\partial x}ds \Big| \mc{G}_t\Big],\quad 0\leq t\leq T, a.s.
\]
where for each $t \in [0,T], (\Gamma_{t,s})_{s \in [t,T]}$ is the unique solution of the following linear SDE
\begin{equation*}
\begin{array}{llll}
 d\Gamma_{t,s} & =&\Gamma_{t,s^-}\Big[\frac{\partial b_{s}}{\partial
x} ds +\frac{\partial
\sigma_{s}}{\partial x} dW_s +\frac{\partial\gamma_{s}}{\partial x}dM_s\Big];\quad \Gamma_{t,t} =1.
\end{array}
\end{equation*}

\subsection{Sufficient Stochastic Maximum Principle}

Now, we are ready to prove the following sufficient maximum
principle for optimal control of an SDE with default of the form \eqref{eq: default}.

\begin{theorem}
\label{thm: suff_max_princ} Let $\hat{u}$ be an admissible performance
strategy with corresponding solution $\hat{X}$ of the SDE
\eqref{eq: default} and the triple adjoint solution $(\hat{p}, \hat{q}, \hat{w})$ to
equation \eqref{eq: adjoint}. Assume

\begin{enumerate}
\item[$(i)$] {The functions } $x \rightarrow g(x)$ and $(x, u) \rightarrow
\mathcal{H}(t, x, u, \hat{p}, \hat{q}, \hat{w})$ are concave a.s. for every $t
\in[0,T]$.

\item[$(ii)$] {For every $v \in\mathcal{V}$},%

\[
 \max_{v\in \mathcal{V}}\mathcal{H}(t, {X}_{t},v,
\hat{p}, \hat{q}, \hat{w})=\mathcal{H}(t, \hat{X}_{t}, \hat{u}_{t},
\hat{p}, \hat{q}, \hat{w}) , \mbox{ }
\mbox{ } dt \times P \mbox{ } a.s.
\]

\end{enumerate}

Then, $\hat{u}$ is an optimal control for the stochastic optimal control
problem with default.
\end{theorem}

Before we move on to the proof, note that the theorem says that if $g$ and the
Hamiltonian are concave, then we may maximize the Hamiltonian instead of the
performance functional in order to find the optimal control of our problem.
This essentially reduces the stochastic optimal control problem to the problem
of solving the SDE \eqref{eq: default} and the adjoint BSDE \eqref{eq: adjoint}. The idea of the proof is to show that
\[
J(u) - J(\hat{u})\leq 0.
\]
From this, the maximum principle follows.

\medskip

\dproof
Fix $\hat{u} \in \mathcal{A}$ with corresponding $\hat{X}_t, \hat{b}_t, \hat{\sigma}_t, \hat{\gamma}_{t}, \hat{p}_t, \hat{q}_t, \hat{w}_t$. Write 
\[
J(u) - J(\hat{u}) = A_1 + A_2,
\]

\noindent where

\[
\begin{array}{llll}
A_1 := E\Big[ \int_0^T \{ h(t, X_t, u_t) - h(t, \hat{X}_t, \hat{u}_t) \}dt\Big], 
A_2 := E\Big[ g(X_T) - g(\hat{X}_T)\Big].
\end{array}
\]
We will use the following notations:'
\small
\[
\mathcal{H}_t:= \mathcal{H}(t, X_t, u_t, \hat{p}_t, \hat{q}_t, \hat{w}_t), \hat{\mathcal{H}}_t :=\mathcal{H}(t, \hat{X}_t, \hat{u}_t, \hat{p}_t, \hat{q}_t, \hat{w}_t),
\]
\[
\frac{\partial\hat{\mathcal{H}}_t}{\partial x}:=\frac{\partial\hat{\mathcal{H}}}{\partial x} (t, \hat{X}_t, \hat{u}_t, \hat{p}_t, \hat{q}_t, \hat{w}_t).
\]
Then, from the definition of the Hamiltonian
\begin{equation}
\label{eq: A_1}
\begin{array}{llll}
A_1 &=& E\Big[ \int_0^T \Big\{ \mathcal{H}_t -  \hat{\mathcal{H}}_t   -  ( b_t- \hat{b}_t)\hat{p}_t -(\sigma_t - \hat{\sigma}_t) \hat {q}_t - \lambda_t(\gamma_{t}- \hat{\gamma}_{t})  \hat{w}_t  \Big\} dt\Big] \\[\medskipamount]
&\leq& E \Big [\int_0^T \Big\{ \frac{\partial \hat{\mathcal{H}}_t}{\partial x} (X_t- \hat{X}_t) +  \frac{\partial \hat{\mathcal{H}}_t}{\partial u} (u_t- \hat{u}_t) - \big( b_t- \hat{b}_t \big)\hat{p}_t  - (\sigma_t - \hat{\sigma}_t)\hat{q}_t  \\[\smallskipamount]
&&- \lambda_t (\gamma_{t}- \hat{\gamma}_{t}) \hat{w}_t \Big\} dt \Big] ,
\end{array}
\end{equation}
\noindent where the second equality follows because the Hamiltonian $\mathcal{H}$ is concave.
Similarly from the concavity of $g$, we have%
\[
\begin{array}{lll}
A_{2} =E\Big[  g(X_{T})-g(\hat{X}_{T})\Big] \leq E \Big[  g^{\prime}(\hat{X}_{T})\left(  X_{T}-\hat
{X}_{T}\right)  \Big]
= E \Big[  \hat{p}_{T}\left(  X_{T}-\hat{X}_{T}\right)
\Big],  \nonumber\\
\end{array}
\]
\noindent where we have used the terminal condition of the BSDE \eqref{eq: adjoint} in the final equality. From It\^{o}'s product rule,
\begin{equation}
\label{eq: A_2'}
d(\hat{p}_{t}(  X_{t}-\hat{X}_{t}))= (X_t - \hat{X}_t) d\hat{p}_t +\hat{p}_t d(X_t - \hat{X}_t) +  d[\hat{p} , (X - \hat{X})]_t. 
\end{equation}
We have 
\[
d[\hat{p},(X-\hat{X})]_t=(\sigma_t - \hat{\sigma}_t)  \hat{q}_tdt+ (\gamma_{t} - \hat{\gamma}_{t})  \hat{w}_tdH_t.
\]
Substituting this into \eqref{eq: A_2'}, integrating between $0$ and $T$ and taking the expectation, we get
\begin{equation}
\label{eq: A_2}
\begin{array}{lll}
E\Big[  \hat{p}_{T}(  X_{T}-\hat{X}_{T})\Big] &=& E\Big[\int_0^T (X_t - \hat{X}_t) (-\frac{\partial\mathcal{\hat{H}}}{\partial x}(t)+\hat{q}_{t}dW_{t}+\hat{w}_{t}dM_{t}) \\[\medskipamount]
&&+ \int_0^T \hat{p}_t \{ \big( b_t- \hat{b}_t\big) dt + (\sigma_t - \hat{\sigma}_t)dW_t+ (\gamma_{t}- \hat{\gamma}_{t}) dM_t \}\\[\medskipamount]
&&
+ \int_0^T\{ \hat{q}_t(\sigma_t - \hat{\sigma}_t)+\lambda_t (\gamma_{t}- \hat{\gamma}_{t})\hat{w}_t\}dt +  \int_0^T (\gamma_{t} - \hat{\gamma}_{t})  \hat{w}_tdM_t\Big].
\end{array}
\end{equation}
Setting
\begin{equation*}
\begin{array}{lll}
d\mathcal{M}_t &=&\Big[(  X_{t}-\hat{X}_{t})\hat{q}_{t}+(\sigma_t - \hat{\sigma}_t)\hat{p}_{t}\Big]dW_t \\[\medskipamount]
&&+\Big[(  X_{t}-\hat{X}_{t})\hat{w}_{t}+(\gamma_{t}- \hat{\gamma}_{t})\hat{p}_{t}+(\gamma_{t}- \hat{\gamma}_{t})\hat{w}_{t}\Big]dM_t.
\end{array}
\end{equation*}
Since $X,\hat{X}, \hat{p} \in S^{2}, \hat{q} \in H^{2}, \hat{w} \in H^{2}(\lambda)$ and the conditions on $\sigma,\hat{\sigma}$ and $\gamma,\hat{\gamma}$, we get that the local martingale $\mathcal{M}$ is a martingale which  has $0$ mean. \\
By combining the expressions for $A_1$ and $A_2$ found in equations \eqref{eq: A_1} and \eqref{eq: A_2} respectively, we find that
\begin{equation}
\label{eq: A1_plus_A2}
A_1 + A_2 \leq E\Big[\int_0^T \frac{\partial \hat{\mathcal{H}}_t} {\partial u} (u_t- \hat{u}_t) dt\Big] \leq 0, 
\end{equation}
\noindent where the final inequality follows from the assumptions. Hence, $J(u) \leq J(\hat{u})$, so since $\hat{u} \in \mathcal{A}$ it is an optimal control.

\fproof

\subsection{Equivalence Maximum Principle}

A problem with the sufficient maximum principle from the previous section is
that the concavity condition is quite strict, and may not hold in
applications. In this section, we derive an alternative maximum principle,
called a necessary maximum principle or equivalence principle for the optimal
control of the SDE with default.

In order to do this we need some additional notation and assumptions:\\
For $u \in \mc{A}$, let $\mc{W}(u)$ denote the set of bounded and $\mathbb{G}$-predictable processes $\beta$ of finite variation such that there exists $\delta=\delta(u) > 0$ satisfying
\begin{equation}
u + y \beta \in \mc{A}, \text{ for all  } y \in [0,\delta].
\end{equation}
Define $\beta:=\mathbf{1}_{[s,T]}(t)\kappa, \text{ for all  } s \in [0,T]$ and $\kappa$ is  bounded and $\mathcal{G}_{s}$-measurable random variable. Assume that for all $u \in \mc{A}$ and for all $\beta \in \mc{W}(u)$ the following derivative process exists and belongs to $L^2([0,T] \times \Omega)$:
\begin{equation}
\label{eq33}
\begin{array}{lll}
x_t:=& \frac{d}{dy} X_t^{u+y\beta}\Big|_{y=0} = \lim_{y \rightarrow 0^+} \frac{X_t^{u + y \beta} - X_t^{u}}{y}. \\
\end{array}
\end{equation}
\begin{remark}
The existence and $L^2$-features of these derivative process is a non-trivial issue, and we do not discuss conditions for this in our paper. We refer to Bensoussan \cite{ben}, Part I, for a study of this issue in a related setting.
\end{remark}
Now, note that from the SDE \eqref{eq: default}, we define the equation of the derivative process
\begin{equation}
\label{der}
\begin{array}{lll}
dx_t=\Big[\frac{\partial b_t}{\partial x}x_t + \frac{\partial b_t}{\partial u}\beta_t \Big] dt+ \Big[\frac{\partial \sigma_t}{\partial x} x_t + \frac{\partial \sigma_t}{\partial u} \beta_t\Big] dW_t + \Big[\frac{\partial \gamma_t}{\partial x} x_{t^{-} }+ \frac{\partial \gamma_t}{\partial u} \beta_t\Big] dM_t; \quad x_0=0.
\end{array}
\end{equation}

\noindent We remark that this derivative process is a linear SDE, then by
assuming that $b$, $\sigma$ and $\gamma$ admit bounded partial derivatives
w.r.t. $x$ and $u$, there is a unique solution $x_t\in
{S}^{2}$ of $\left(  \ref{der}\right). $
The following equivalence principle says that for a control to be a critical point for the performance functional $J$ is equivalent to the critical point of the Hamiltonian.

\begin{theorem}
\label{thm: equivalence_principle} The following two statements are equivalent:

\begin{enumerate}
\item[$(i)$] {
\begin{equation}
\label{eq: i}\frac{d J(u+y\beta)}{dy}\Big|_{y=0} = 0.
\end{equation}
}

\item[$(ii)$] {
\[
\frac{\partial\mathcal{H}_t}{\partial u} = 0.
\]
}
\end{enumerate}
\end{theorem}

\dproof
Note that
\[
\frac{d J(u+y\beta)}{dy}\Big|_{y=0} = \frac{d}{dy} E\Big[\int_0^T h(t, X_t^{u + y\beta}, u_t+y_t \beta_t)dt + g(X_T^{u+y\beta})\Big]\Big|_{y=0}  .
\]
Define $I_1 := \frac{d}{dy}E\Big[\int_0^T h(t, X_t^{u + y\beta}, u_t+y\beta_t)dt\Big]\Big|_{y=0}$ and $I_2 :=\frac{d}{dy} E\big[g(X_T^{u+y\beta})\big]\Big|_{y=0} $.
Since the coefficients have uniformly bounded derivatives, it follows from the dominated convergence theorem that the equality follows
\[
I_1 = E\Big[\int_0^T \Big\{ \frac{\partial h_t}{\partial x}x_t + \frac{\partial h_t}{\partial u}\beta_t \Big\} dt\Big].
\]
Also,
\[
\begin{array}{lll}
I_2 &=& E\big[g'(X_T^{u+y\beta}) x_T\Big]= E[p_Tx_T],
\end{array}
\]
\noindent where the first equality follows by changing the order of differentiation and integration (again, using the dominated convergence theorem), the second equality follows from the adjoint equation \eqref{eq: adjoint}. 
So, by the previous expression for $dx_t$ \eqref{der} and $d[p,x]_t$, as well as the expression for $dp_t$ from the BSDE \eqref{eq: adjoint}, we obtain
\[
\begin{array}{llll}
I_2 &=& E\Big[\int_0^T p_t ( \frac{\partial b_t}{\partial x} x_t + \frac{\partial b_t}{\partial u} \beta_t)dt - \int_0^T x_t \frac{\partial \mathcal{H}_t}{\partial x}dt + \int_0^T q_t (\frac{\partial \sigma_t}{\partial x} x_t + \frac{\partial \sigma_t}{\partial u}\beta_t)dt\\[\smallskipamount]
&& + \int_0^T \lambda_t w_t (\frac{\partial \gamma_t}{\partial x} x_t + \frac{\partial \gamma_t}{\partial u} \beta_t)dt\Big].
\end{array}
\]
By collecting the $\beta_t$- and $x_t$-terms and using the definition of the Hamiltonian to cancel all $x_t$-terms against $x_t \frac{\partial \mathcal{H}_t}{\partial x}$, we find that
\[
\begin{array}{llll}
I_1 + I_2 &=& E\Big[\int_0^T \beta_t (\frac{\partial h_t}{\partial u} + p_t\frac{\partial b_t}{\partial u} + q_t \frac{\partial \sigma_t}{\partial u} + w_t \lambda_t \frac{\partial \gamma_t}{\partial u}) dt\Big] \\[\smallskipamount]
&=& E\Big[\int_0^T \beta_t \frac{\partial \mathcal{H}_t}{\partial u}dt\Big],  \mbox{  for all } \beta \in \mc{W}(u).
\end{array}
\]
In particular, if we apply this to
\[
\beta=\mathbf{1}_{[s,T]}(t)\kappa,
\]
where $\kappa$ is bounded and $\mathcal{G}_{s}$-measurable, we get
\[
0\geq\mathbb{E}\left[  \int_{s}^{T}\frac{\partial\mathcal{H}_t}{\partial
u}\kappa dt\right]  .
\]
Since this holds for all such $\kappa$ (positive or negative) and all
$s\in[0,T]$, we conclude that
\[
0 = \frac{\partial\mathcal{H}_t}{\partial u},\quad\text{ for a.a. }t,
\]
\noindent and hence the theorem follows.
\fproof
\section{Application: Log-Utility Maximisation with Default}

In this section, we illustrate the stochastic maximum principles Theorem \ref{thm: suff_max_princ} and Theorem \ref{thm: equivalence_principle} by applying them to a logarithmic utility maximisation problem. As pointed out by the referee, it is also possible to solve this problem directly by using the formula of $X_t$ and integrating by parts.\\
Consider the cash flow process with default
\begin{equation}
\label{eq: SDE_ex}
dX_{t}=X_{t-}\Big[(\alpha_{t}-c_{t})dt+\rho_{t}dW_{t}+\mu_{t}dM_{t}\Big];\quad X_{0}>0,
\end{equation}
where the coefficients $\alpha,\rho,\mu$ are bounded, $\mathbb{R}$-valued $\mathbb{G}$-predictable
processes and we assume that $\mu_{t}\geq-1$ for all $t \in [0,T]$ a.s. From the so-called, Dol\'{e}ans-Dade formula, we can write the linear SDE \eqref{eq: SDE_ex} explicitly, as follows
\small
\begin{align}
    X_t=X_0 exp\Big(\int_{0}^t \Big\{\alpha _s-c_s-\frac{1}{2}\rho^2_s\Big\}ds+\int_{0}^t\rho_sdW_s\Big)exp\Big(-\int_{0}^t \mu_s \lambda_s ds \Big) \Big(1+\mu_\tau\mathbf{1}_{\tau\leq t}\Big).
\end{align}
Since $\mu_{\tau}\geq-1$ and $X_{0}>0$ imply that $X_{t}>0$ a.s. for each $t\in [0,T]$.
Also, note that in the SDE \eqref{eq: SDE_ex}, the control $c_t\geq0$ corresponds to a consumption process because of its negative impact on the cash flow process $X_t$. The default term $\mu_{t}dM_{t}$ implies that the wealth process will grow w.r.t. $\mu_t$ until the default time $\tau$. From that point on, $\mu_t$ has no impact on the cash flow. This may correspond to investing in a defaultable firm.
The performance function we want to maximize is%

\[
J(c)=E\Big[%
{\textstyle\int_{0}^{T}}
U^{1}(X_{t},c_{t})dt+\theta U^{2}( X_{T})\Big],
\]

\noindent where $U^{1},U^{2}$ are some given deterministic utility functions and $\theta
:=\theta(\omega)>0$ is a $\mathcal{G}_{T}$-measurable, square integrable random
variable which expresses the importance of the terminal value. To be able to
find explicit solutions for our optimal control, we consider logarithmic utilities. Hence, the performance function is
\[
J(c)=E\Big[%
{\textstyle\int_{0}^{T}}
\log(X_{t}c_{t})dt+\theta\log( X_{T})\Big].
\]
The corresponding Hamiltonian functional, see \eqref{eq: Hamiltonian}, is
\[
\mathcal{H}(t,x,c,p,q,w)=\log(xc)+x(\alpha-c)p+x\rho
q+\lambda_t x\mu w.
\]
The adjoint BSDE, see \eqref{eq: adjoint}, has the form%
\begin{equation}
\label{eq: BSDE_example}
dp_{t}  =-\frac{\partial\mathcal{H}_t}{\partial x}dt+q_{t}dW_{t}+w_{t}%
dM_{t};\quad
p_{T}  =\frac{\theta}{X_{T}},
\end{equation}
such that%
\begin{equation}
\label{eq: partial_H}
\frac{\partial\mathcal{H}_t} {\partial x}=\frac{1}{X_{t}}+(\alpha_{t}-c
_{t})p_{t}+\rho_{t}q_{t}+\lambda_{t}\mu_{t}w_{t}.
\end{equation}
Using the first order necessary condition of Theorem \ref{thm: equivalence_principle}, we obtain%
\[
\frac{\partial\mathcal{H}_t} {\partial c}=\-X_{t}p_{t}+\frac{1}{c_{t}}=0.
\]
Consequently,%
\begin{equation*}
\hat{ c_{t}}=\frac{1}{X_{t}p_{t}}.\label{pi}%
\end{equation*}
To derive an explicit expression for the optimal control, we do the following
computations:

Note that by the It{\^o} product rule,

\[
d(p_t X_t) = p_t d X_t + X_t dp_t + d[p,X]_t.
\]

Hence,

\begin{equation}
\label{eq: product}
\begin{array}{lll}
p_T X_T - p_t X_t &=& \int_t^T X_s\Big(p_s(\alpha_s-c_s)-\frac{\partial\mathcal{H}_s} {\partial x}+w_s\mu_s\lambda_s\Big)ds+X_s\Big(p_s\rho_s+q_s\Big)dW_s \\[\medskipamount]
&&+X_{s-}\Big(p_{s-}\mu_s+w_s(1+\mu_s)\Big)dM_s.
\end{array}
\end{equation}

By taking the conditional expectation w.r.t $\mc{G}_t$ on both sides of equation \eqref{eq: product}, we see that

\[
\begin{array}{lllll}
p_t X_t &=& E\Big[p_T X_T \Big| \mc{G}_t\Big] + E\Big[\int_t^T  ds \Big| \mc{G}_t\Big] \\[\smallskipamount]
&=& E\Big[\theta+T-t \Big| \mathcal{G}_{t}\Big],
\end{array}
\]
\noindent where the second to last equality follows by inserting expression for $p_T$ from the adjoint BSDE \eqref{eq: BSDE_example}.


Hence, an explicit expression for the stochastic optimal consumption is
\begin{equation*}
\hat{ c_{t}}=\frac{1}{E\Big[\theta+T-t \Big| \mathcal{G}_{t}\Big]} = \frac{1}{E[\theta | \mathcal{G}_t] + T - t}.\label{pi_explicit}%
\end{equation*}

\thanks{{\bf Acknowledgements}. We are grateful to an anonymous reviewer for valuable comments and feedback.
N. Agram and K. Dahl are gratefully acknowledge the financial support provided by the Swedish Research Council grant (2020-04697) and the Norwegian Research Council grant (299897), respectively.}


\begin{thebibliography}{99}                                                                                               

\bibitem {aj}Aksamit, A., \& Jeanblanc, M. (2017). Enlargement of filtration
with finance in view. Springer.
\bibitem {ben}Bensoussan, A. (1982). Lectures on stochastic control. In Nonlinear filtering and stochastic control (pp. 1-62). Springer, Berlin, Heidelberg.


\bibitem{BouchardPham}
Bouchard, B., \& Pham, H. (2004). Wealth-path dependent utility maximization in incomplete markets. Finance and Stochastics, 8(4), 579-603.

\bibitem{ChoulliYansori}
Choulli, T., \& Yansori, S. (2018). Deflators and log-optimal portfolios under random horizon: Explicit description and optimization. arXiv preprint arXiv:1803.10128.

\bibitem{CohenElliott}
Cohen, S. N., \& Elliott, R. J. (2015). Stochastic calculus and applications (Vol. 2). New York: Birkh\"{a}user.

\bibitem{DumitrescuEtAl}
Dumitrescu, R., Grigorova, M., Quenez, M. C., \& Sulem, A. (2016, August). BSDEs with default jump. In The Abel Symposium (pp. 233-263). Springer, Cham.

\bibitem{EyjolfssonTjostheim}
Eyjolfsson, H., \& Tj{\o}stheim, D. (2018). \emph{Self-exciting jump processes with applications to energy markets}, Ann Inst Stat Math, 70: 373-393.

\bibitem{Hawkes1}
Hawkes, A. G. (1972). Spectra of some self-exciting and mutually exciting processes, Biometrika, 58, 83-90.

\bibitem{Hawkes2}
Hawkes, A. G., \& Oakes, D. (1974). A cluster process representation of a self-exciting process, Journal of Applied Probability, 11, 493-503.

\bibitem{JeanblancEtAl}
Jeanblanc, M., Mastrolia, T., Possama\"{i}, D., \& R\'{e}veillac, A. (2015). Utility maximization with random horizon: a BSDE approach. International Journal of Theoretical and Applied Finance, 18(07), 1550045.

\bibitem {kha}Kharroubi, I., \& Lim, T. (2014). Progressive enlargement of
filtrations and backward stochastic differential equations with jumps. Journal
of Theoretical Probability, 27(3), 683-724.

\bibitem{LimQuenez}
Lim, T., \& Quenez, M. C. (2011). Exponential utility maximization in an incomplete market with defaults. Electronic Journal of Probability, 16, 1434-1464.

\bibitem{LimQuenez2}
Lim, T., \& Quenez, M. C. (2015). Portfolio optimization in a default model under full/partial information. Probability in the Engineering and Informational Sciences, 29(4), 565-587, 

\bibitem {P}Pham, H. (2010). Stochastic control under progressive enlargement of filtrations and applications to multiple defaults risk management. Stochastic Processes and their applications, 120(9), 1795-1820.


\end{thebibliography}
\end{document}